\documentstyle{article}
\title{LATTICE COVERINGS AND GAUSSIAN MEASURES OF $n$-DIMENSIONAL CONVEX BODIES}
\author{Wojciech Banaszczyk\thanks{Part of this research was
done while
this author was visiting Case Western Reserve University under a cooperation
grant from KBN (Poland) and NSF (U.S.A.)} (\L\'od\'z)
 and Stanislaw J. Szarek\thanks{Supported in part by the
National Science Foundation.} (Cleveland)}
\date{}
\begin{document}
\maketitle

\begin{abstract}
Let $\| \cdot \|$ be the euclidean norm on ${\bf R}^n$ and $\gamma_n$ the
(standard)
Gaussian measure on ${\bf R}^n$ with density $(2 \pi )^{-n/2} e^{- \| x\|^2
/2}$.  Let
$\vartheta$  ($ \simeq 1.3489795$) be defined by $\gamma_1 ([ -  \vartheta
/2,  \vartheta /2]) = 1/2$ and let $L$ be a
lattice in ${\bf R}^n$ generated by vectors of norm $\leq  \vartheta$.
Then, for any closed
convex set $V$ in ${\bf R}^n$ with $\gamma_n (V) \geq \frac{1}{2}$ and for any
$a \in {\bf R}^n$, $(a +L) \cap V \neq \phi$.  The above statement can be
viewed as a
``nonsymmetric'' version of Minkowski Theorem.
\end{abstract}

\bigskip Let $U$, $V$ be a pair of convex sets in ${\bf R}^n$ containing the
origin in the interior.
Let us define $\beta (U,V)$ as the smallest $r > 0$ satisfying the
following condition: to each
sequence $u_1 , \ldots , u_n \in U$ there correspond signs $\varepsilon_1 ,
\ldots ,
\varepsilon_n = \pm 1$ such that $\varepsilon_1 u_1 + \cdots +
\varepsilon_n u_n \in rV$.  Upper
and lower bounds for $\beta (U,V)$ for various sets $U$ and $V$ (usually
centrally symmetric) were
investigated by several authors. We will mention some of their results once the
appropriate notation is introduced, see also references in [3].

Let $L$ be a lattice in ${\bf R}^n$, i.e. an additive subgroup of ${\bf
R}^n$ generated by $n$
linearly independent vectors.  The quantities (again,  usually defined for
centrally
symmetric sets)
$$\lambda_n (L,U) = \min \{ r >0: \mbox{ dim span } (L \cap rU) =n \} ,$$

$$\mu (L,V) = \min \{ r > 0: L + rV = {\bf R}^n \}$$
are called the $n$th minimum and the covering radius of $L$ with respect to
$U$ and $V$,
respectively;  sometimes  $\mu (L,V)$  is called "the $n$th covering
minimum" and denoted  $\mu_n (L,V)$.  Let us define

\noindent \hrulefill

\noindent AMS Subject Classification 11H06, 11H31, 52C07, 52C17
\newpage

$$\alpha (U,V) = \sup_L \frac{\mu (L,V)}{\lambda_n (L,U)}$$
where the supremum is taken over all lattices $L$ in ${\bf R}^n$.  A
standard elementary
argument shows that $\alpha (U,V) \leq \beta (U,V)$  (see e.g. Lemma 4 in
[3]).

By $B_n$ we shall denote the closed euclidean unit ball in ${\bf R}^n$.
Let $E$ be an
$n$-dimensional ellipsoid in ${\bf R}^n$ with centre at zero and principal
semiaxes
$\alpha_1 , \ldots , \alpha_n$.  The result of [4], that closed connected
additive subgroups of
nuclear locally convex spaces are linear subspaces, was essentially based
on the fact that
$$\alpha (B_n, E) = \frac{1}{2} (\alpha_1^2 + \cdots + \alpha_n^2 )^{1/2}.$$
Then it was proved in [2] that
$$\beta (B_n, E) = (\alpha_1^2 + \cdots + \alpha_n^2)^{1/2}.$$

Let $K_n$ be the unit cube in ${\bf R}^n$.  Consider the rectangular
parallelepiped
$$P = \{ (x_1 , \ldots , x_n) \in {\bf R}^n :
|x_k| \leq \alpha_k \mbox{ for } k =1, \ldots ,n \}$$
where $\alpha_1 , \ldots , \alpha_n >0$.  This paper was motivated by an
attempt to give
possibly best upper bounds for $\alpha (B_n ,P)$ and $\beta (B_n ,P)$ as
functions of
$\alpha_1 , \ldots , \alpha_n$ (for $\beta (K_n,P)$, see [5] and [8] where
it was, in
particular, proved that $\beta (K_n , K_n) =O (\sqrt{n} )$ as $n
\rightarrow \infty$;  see also  [1]).
In particular, we were interested in the so-called Koml\'{o}s conjecture
which asserts
that $\beta (B_n , K_n)$ remains bounded as $n \rightarrow \infty$.

Let us denote by $\gamma_n$ the (standard) Gaussian measure on ${\bf R}^n$
with density
$(2 \pi )^{-n/2} e^{- \| x \|^2 /2}$, where $\| x \|$ is the euclidean norm
of $x$.  Let
$\vartheta$  ($ \simeq 1.3489795$)
be the positive number given by $\gamma_1 ([- \vartheta /2,
\vartheta /2]) = \frac{1}{2}$, i.e.

$$\int_0^{\vartheta /2} e^{-t^2 /2} dt = \frac{\sqrt{2 \pi}}{4} .$$
By a $\vartheta$-coset in ${\bf R}^n$ we shall mean a coset modulo a lattice $L$
generated by vectors of Euclidean norm $\leq  \vartheta$, i.e. satisfying
$\lambda_n (L, B_n ) \leq \vartheta$.  The aim of this paper is to prove
the following
fact.

\medskip\noindent {\bf Theorem.} If $V$ is a closed convex set in ${\bf R}^n$
with $\gamma_n (V) \geq 1/2$, then $V$ intersects every $ \vartheta$-coset.

\medskip\noindent {\bf Corollary.} If $V$ is as in the Theorem, then
$\alpha (B_n, V) \leq  \vartheta^{-1}$.  In particular $\alpha (B_n , K_n ) = O
\sqrt{\log n})$
as $n \rightarrow \infty$

\medskip We point out that, in full generality, the Theorem is sharp and that,
similarly,  the
first part of the Corrollary can not be significantly improved.  However, it is
conceivable that $\alpha (B_n , \cdot )$ may be replaced by $\beta (B_n ,
\cdot )$ in the
Corollary; see the Conjecture at the end of this paper.

\medskip For the proof we need the following.

\medskip\noindent {\bf Lemma.} If $V$ is a closed convex set in ${\bf R}^n$ with
$\gamma_n (V) \geq \frac{1}{2}$ and $M$ is a linear subspace of ${\bf R}^n$
of dimension $m$,
then $\gamma_m (V \cap M) \geq \frac{1}{2}$.

\medskip\noindent {\bf Remark 1.} An analysis of the proof shows that unless
$V$ is a half space,  or an infinite cylinder orthogonal to $M$,
the inequality in the assertion of the Lemma is strict.

\medskip We need some preparation for the proofs of the Lemma and of the
Theorem.
For a convex
set $V$ in ${\bf R}^n$ and $x \in {\bf R}$ denote

\begin{equation}
V_x = \{ (x_1 , \ldots , x_{n-1} ) \in {\bf R}^{n-1} : (x_1 , \ldots ,
x_{n-1} ,x) \in V \}
\end{equation}
Recall now an inequality of Ehrhard (see [6], Thm. 3.2).  If $A$, $B$
are non-empty convex
Borel subsets of ${\bf R}^n$  and  $0 \leq \lambda \leq 1$, then
\begin{equation}
\Phi^{-1} ( \gamma_n ( \lambda A + (1- \lambda )B)) \geq \lambda \Phi^{-1}
( \gamma_n (A)) +
(1- \lambda ) \Phi^{-1} ( \gamma_n (B))
\end{equation}
where

$$\Phi (x) = \frac{1}{\sqrt{2 \pi}} \int_{- \infty}^x e^{-y^2 /2} dy,
\qquad x \in {\bf R}.$$
is the (standard) Gaussian cumulative distribution function.  It follows in
particular that
$g(x) = \Phi^{-1} ( \gamma_{n-1} (V_x))$ is a concave function of $x$ on the
interval $I = \{ x: \gamma_{n-1} (V_x) > 0 \}$.  Consequently,
\begin{equation}
W = \{ (x,y) \in {\bf R}^2 : x \in I \mbox{ and } y \leq g(x) \}
\end{equation}
is a closed convex subset of ${\bf R}^2$.  Note that $\gamma_1 (W_x) =
\gamma_1 ((
- \infty ,g(x)]) = \gamma_{n-1} (V_x)$ for $x \in {\bf R}$, where $W_x$ is
defined
analogously to $V_{x}$;  in particular $\gamma_n (V) = \gamma_2 (W)$.

\medskip\noindent {\bf Proof of the Lemma.}  Clearly it is enough to
consider the case
$m =n -1$ and (by the rotationary invariance of the Gaussian measure)
$M = \{ (x_1 , \ldots , x_n) : x_n =0 \}$.  For $V$ with $\gamma_n (V) \geq
\frac{1}{2}$ we construct $W \subset {\bf R}^2$ as above, the assertion of the
Lemma is then equivalent to $\gamma_1 (W_0) \geq \frac{1}{2}$
or $(0,0) \not\in W$.  To conclude the argument it remains to note that $(0,0)
\not\in W$, together with
$W$ being closed and convex, would imply $\frac{1}{2} > \gamma_2 (W) = \gamma_n
(V)$, a contradiction.

\medskip\noindent {\bf Remark 2.} For the proof of the Theorem we use the Lemma
with $n =2$ and $m =1$,  a special case that can be proved without
appealing to the
Ehrhard's inequality (2).  However, the proof of the Theorem itself does
use Ehrhard's
inequality.

\medskip\noindent {\bf Proof of the Theorem.} We use induction on $n$.  For $n
=1$, the Theorem is rather trivial.
So, suppose that for a certain $n \geq 2$ the Theorem is true for all
dimensions strictly less than $n$ .  Take an arbitrary
$\vartheta$-coset
$H$ in ${\bf R}^n$ and a convex set $V$ in ${\bf R}^n$ disjoint with $H$.
We are to prove
that $\gamma_n (V) < \frac{1}{2}$.

Fix some $u \in H$ and consider the lattice $L =H -u$.  By assumption, we have
$\lambda_n (L,B_n) \leq \vartheta$.  Choose $a_1 , \ldots , a_n \in L \cap
 \vartheta B_n$
generating $L$ and let $M$ be the linear span of $a_1 , \ldots , a_{n-1}$.
As before,
we may assume that $M = \{ (x_1 , \ldots , x_n): x_n =0 \}$.  Let $H'$ be
the orthogonal projection
of $H$ onto the $n$th coordinate axis of ${\bf R}^n$ (i.e., onto the orthogonal
complement of $M)$.  Clearly $H'$ is a $\vartheta$-coset.  Additionally, if
$x \in H'$, then, by
our inductive hypothesis, $\gamma_{n-1} (V_x) < \frac{1}{2}$ and so $(x,0)
\not\in
W$ $(V_x ,W$ have
the samemeaning here as in (1) and (3)).  The case $n =1$ of the Theorem yields
now that
$\gamma_1 (W \cap \{ (x,0) : x \in {\bf R} \} ) < \frac{1}{2}$ and the
Lemma implies then that
$\frac{1}{2} > \gamma_2 (W) = \gamma_n (V)$, as required.

\medskip\noindent {\bf Conjecture.} There exists some function $f$ on
$(0,1)$ such that for each
symmetric convex set $V$ in ${\bf R}^n$ one has $\beta (B_n ,V) \leq f(
\gamma_n (V))$.

\medskip\noindent {\bf Remark 3.} Let $T$ be a bounded linear operator from a
Hilbert space $H$ to a
Banach space $X$.  We say that $T$ is tight if the image of every connected
additive subgroup
of $H$ is dense in its linear span in $X$.  If $X$ is a Hilbert space, then
$T$ is tight if and only
if it is a Hilbert-Schmidt operator; sufficiency was proved in [4], the
proof of necessity
can easily be obtained by standard methods.  The argument of [4] together
with the theorem
proved above imply that $\ell$-operators are tight (for the definition of
$\ell$-operators, see [7],
p. 38).  An interesting problem, closely connected with the Koml\'{o}s
conjecture, is to
describe tight diagonal operators from $l_2$ to $c_0$.

\end{document}